\documentclass{amsart}[12pt]
\usepackage{latexsym}
\usepackage{amsmath}
\usepackage{amssymb}
\usepackage{amsfonts}
\usepackage{amsthm}

\newcommand {\Z}  {\ensuremath{\mathbb{Z}}}

\newcommand{\raisedLabel}[1]{\leavevmode \raise 12pt\hbox to 0pt{\tiny\text{#1}\hss}\hskip-2pt}

\def\newNumber#1{%
\refstepcounter{thm}\label{#1}}

\newcommand{\namedRef}[1]{%
\csname env:#1\endcsname\ \ref{#1}%
}

\def\envLRT#1#2{\expandafter\gdef\csname env:#1\endcsname{#2}}

\makeatletter
\newcommand{\namedLabel}[1]{\expandafter\xdef\csname env:#1\endcsname{\envName}%
\immediate\write\@auxout{\string\envLRT{#1}{\envName}}%
\label{#1}\raisedLabel{#1}}
\newcommand{\NamedLabel}[1]{\expandafter\xdef\csname env:#1\endcsname{\envName}%
\immediate\write\@auxout{\string\envLRT{#1}{\envName}}%
\label{#1}}
\makeatother

\newcommand{\newType}[2][]{%
\def\LRTxx{#1}\ifx\LRTxx\empty\else\theoremstyle{#1}\fi%
\newtheorem{#2}[thm]{#2\gdef\envName{#2}}}
\newtheorem{thm}{Theorem}[section]
\newtheorem{Theorem}[thm]{Theorem\gdef\envName{Theorem}}
\newType{Lemma}
\newType{Corollary}
\newType{Proposition}
\newType{Addendum}
\newType{Result}
\newenvironment{Named Theorem}[1][xxx]{%
\gdef\envName{#1}\par\noindent{\bf #1.}\bgroup\sl }{\egroup}

\theoremstyle{definition}

\newtheorem{Example}[thm]{Example\gdef\envName{Example}}
\newtheorem{Remark}[thm]{Remark\gdef\envName{Remark}}

\newType{Notation}
\newType{Examples}

\newif\ifShowMathCheck
\ShowMathCheckfalse

\newtheorem{realMathCheck}[thm]{Check}
\newbox\hidebox
\newenvironment{hiddenMathCheck}{\setbox\hidebox=\vbox\bgroup}
{\egroup\setbox\hidebox=\vbox{}}

{\ifShowMathCheck\relax\begin{realMathCheck}\else\begin{hiddenMathCheck}\fi}
{\ifShowMathCheck\relax\end{realMathCheck}\else\end{hiddenMathCheck}\fi}

\newcommand{\cy}[1]{\Z/{#1}\Z}

\newcommand{\longequal}[1]{%
\vbox{\hsize=#1\hbox to#1{\vrule width #1 height .5pt depth 0pt}\vskip 3pt
\hbox to#1{\vrule width #1 height .5pt depth 0pt}}}
\newcommand{\adjustedrightlabeledarrow}[5][]{%
\hbox to #4{\hskip#5$\rightlabeledarrow[#1]{#2}{#3}$\hss}}
\newcommand{\adjustedlongequal}[3]{%
\hbox to #2{\hskip#3\longequal{#1}}}

\newcommand{\rightlabeledarrow}[3][]{%
\def\LRTxx{#1}\ifx\LRTxx\empty\relax%
\setbox0=\hbox{$\mathop{\expandafter\hbox%
{\rightarrowfill}}\limits_{\hbox{\ $\scriptstyle#3$\ }}^{\hbox{\ $\scriptstyle#2$\ }}$}%
\dimen0=\wd0\advance\dimen0 by 4pt%
\edef\LRTxx{\the\dimen0}\else%
\setbox0=\hbox to \LRTxx{$\mathop{\expandafter\hbox%
{\rightarrowfill}}\limits_{\hbox{$\scriptstyle#3$}}^{\hbox{$\scriptstyle#2$}}$}%
\dimen0=\wd0\advance\dimen0 by 4pt%
\fi%
\setbox0=\hbox{$\mathop{\expandafter\hbox%
{\rightarrowfill}}\limits_{\hbox{$\scriptstyle#3$}}^{\hbox{$\scriptstyle#2$}}$}%
\ifdim\dimen0<\wd0 \relax\dimen0=\wd0 \advance\dimen0 by 4pt\fi%
\ \hbox to\dimen0{\hfill\hbox to 0pt{\hss$\mathop{\hbox to\LRTxx {\rightarrowfill}}%
\limits_{\hbox to 0pt{\hss{$\scriptstyle#3$}\hss}}^{\hbox to 0pt{\hss{$\scriptstyle#2$}\hss}}$\hss}\hfill}\ %
}

\newcommand{\leftlabeledarrow}[3][]{%
\def\LRTxx{#1}\ifx\LRTxx\empty\relax%
\setbox0=\hbox{$\mathop{\expandafter\hbox%
{\leftarrowfill}}\limits_{\hbox{\ $\scriptstyle#3$\ }}^{\hbox{\ $\scriptstyle#2$\ }}$}%
\dimen0=\wd0\advance\dimen0 by 4pt%
\edef\LRTxx{\the\dimen0}\else%
\setbox0=\hbox to \LRTxx{$\mathop{\expandafter\hbox%
{\leftarrowfill}}\limits_{\hbox{$\scriptstyle#3$}}^{\hbox{$\scriptstyle#2$}}$}%
\dimen0=\wd0\advance\dimen0 by 4pt%
\fi%
\setbox0=\hbox{$\mathop{\expandafter\hbox%
{\leftarrowfill}}\limits_{\hbox{$\scriptstyle#3$}}^{\hbox{$\scriptstyle#2$}}$}%
\ifdim\dimen0<\wd0 \relax\dimen0=\wd0 \advance\dimen0 by 4pt\fi%
\ \hbox to\dimen0{\hfill\hbox to 0pt{\hss$\mathop{\hbox to\LRTxx {\leftarrowfill}}%
\limits_{\hbox to 0pt{\hss{$\scriptstyle#3$}\hss}}^{\hbox to 0pt{\hss{$\scriptstyle#2$}\hss}}$\hss}\hfill}\ %
}

\newcommand{\downlabeledarrow}[3][]{%
\raise2pt\hbox to 0pt{\hss$\scriptstyle #2$}%
\def\LRTxx{#1}\ifx\LRTxx\empty\relax\downarrow%
\else #1\downarrow\fi\raise2pt\hbox to 0pt{$\scriptstyle #3$\hss}%
}
\newcommand{\downlabeledarrowadjusted}[5][]{%
\raise#4\hbox to 0pt{\hss$\scriptstyle #2$}%
\def\LRTxx{#1}\ifx\LRTxx\empty\relax\downarrow%
\else #1\downarrow\fi\raise#5\hbox to 0pt{$\scriptstyle #3$\hss}%
}

\newcommand{\uplabeledarrow}[3][]{%
\raise2pt\hbox to 0pt{\hss$\scriptstyle #2$}%
\def\LRTxx{#1}\ifx\LRTxx\empty\relax\uparrow%
\else #1\uparrow\fi\raise2pt\hbox to 0pt{$\scriptstyle #3$\hss}%
}

\newbox\nbox
\gdef\Bar#1{\setbox\nbox = \hbox{$#1$}
\kern .08\wd\nbox{\overline{\hbox to .84\wd\nbox{\vphantom {$#1$}\hss}}}
\kern .08\wd\nbox\kern -\wd\nbox\box\nbox}

\newtoks\LRTta \newtoks\LRTtb
\long\def\LRTleftappenditem#1\to#2{\LRTta={\\{#1}}\LRTtb=\expandafter{#2}%
\edef#2{\the\LRTta\the\LRTtb}} 
\long\def\LRTrightappenditem#1\to#2{\LRTta={\\{#1}}\LRTtb=\expandafter{#2}%
\edef#2{\the\LRTtb\the\LRTta}} 

\gdef\BibList{}
\def\localCite{%
\gdef\localCite##1##2{\LRTrightappenditem{{##1}{##2}}\to\BibList}
\def\\##1{\LRTrr##1@}%
\def\LRTrr##1##2@{\def\gg{##2}\ifx\gg\LRTG\relax\gdef\LRTreturn{##1}\fi}
\let\BiB=\bib
\def\MR##1##2##3@{\gdef\xx{##1##2}}
\renewcommand{\bib}[3]{%
\def\LRTG{##1}\BibList\ifx\LRTreturn\empty\relax\else%
\expandafter\BiB\expandafter{\LRTreturn}{##2}{##3}\fi%
\gdef\LRTreturn{}%
\MR##1@\def\yy{MR}\ifx\xx\yy\else\BiB{##1}{##2}{##3}\fi}
\localCite}

\usepackage{tikz} 
\usetikzlibrary{arrows} 
\usetikzlibrary{positioning,calc,through,decorations.pathmorphing,petri,backgrounds,fit} 

\def\eshow#1{\expandafter\show\csname #1\endcsname}

\def\cs#1#2{\expandafter\def\csname #1\endcsname{#2}}

\def\addCoor#1#2{%
\pgfmathparse{\lastSum+#2}\edef\lastSum{\pgfmathresult}%
\expandafter\edef\csname #1\endcsname{\pgfmathresult}%
}

\newcommand{\deltaX}[1]{\def\deltaXlist{#1}}
\newcommand{\deltaY}[1]{\def\deltaYlist{#1}}

\newcommand{\setCoor}[3][]{%
\gdef\lastSum{#2}%
\def\\##1##2{\addCoor{#1x##1}{##2}}\deltaXlist%
\gdef\lastSum{#3}%
\def\\##1##2{\addCoor{#1y##1}{##2}}\deltaYlist%
}

\let\\=\relax
\gdef\arrowData{[->]}

\newcommand{\DN}[2][]{\def\\##1##2##3{%
\draw node (#1##1##2) at (\csname #1x##1\endcsname, \csname #1y##2\endcsname) {$##3$};%
}#2}
\newcommand{\DhA}[2][]{\def\\##1##2##3{%
\draw \arrowData (#1##1) to node [midway, above] {$\scriptstyle##3$} (#1##2);%
}#2}

\newcommand{\DvA}[2][]{\def\\##1##2##3{%
\draw \arrowData (#1##1) to node [midway, right] {$\scriptstyle##3$} (#1##2);%
}#2}

\def\leqRef#1#2#3#4{%
\draw node at ($(#1) + (0,#2)$) {\hbox to 0pt{\hss(\ref{#4})\hskip #3pt\null}};%
}
\def\leqNO#1#2#3#4{%
\draw node at ($(#1) + (0,#2)$) {\hbox to 0pt{\hss#4\hskip #3pt\null}};%
}

\newtoks\ta\newtoks\tb
\long\def\rightappenditem#1\to#2{%
\ta={\\#1}\tb=\expandafter{#2}%
\xdef#2{\the\tb\the\ta}}

\gdef\LRTarrowKill{.}
\gdef\tikizSetup{\begin{tikzpicture}}

\newcount\ladderType
\newcount\rungs
\newcount\ladderCount
\newcount\tikzLRTdepthCount
\global\tikzLRTdepthCount=0
\newcount\nextRungLRT

\def\dLRTX#1,#2\\#3\\{%
\global\advance \ladderCount by 1 %
\ifnum\ladderType=1 %
\expandafter\rightappenditem\expandafter{\number\ladderCount}{#1}\to\deltaXlist%
\else%
\expandafter\rightappenditem\expandafter{\number\ladderCount}{-#1}\to\deltaYlist\fi%
\def\xdLRTXx{#2}\ifx\xdLRTXx\empty %
\global\rungs\ladderCount %
\expandafter\setCoor\LRToption \global\ladderCount=0 \let\next= \dLRTRungs\def\xdLRTXx{#3\\}%
\else%
\let\next=\dLRTX\def\xdLRTXx{#2\\#3\\}\fi\expandafter\next\xdLRTXx%
}

\def\dLRTRungs#1#2#3#4\\#5\\{%
\global\advance\ladderCount by 1 %
\def\xLRTRungsx{#2}
\ifnum\ladderType=1 %
\DN[\LRToptionName]{\\{\number\ladderCount}1{#1}\\{\number\ladderCount}2{#3}}%
\ifx\xLRTRungsx\LRTarrowKill \else\DvA[\LRToptionName]{\\{\number\ladderCount1}{\number\ladderCount2}{#2}}\fi%
\else%
\DN[\LRToptionName]{\\1{\number\ladderCount}{#1}\\2{\number\ladderCount}{#3}}%
\ifx\xLRTRungsx\LRTarrowKill \else\DhA[\LRToptionName]{\\{1\number\ladderCount}{2\number\ladderCount}{#2}}\fi%
\fi%
\ifnum\ladderCount=\rungs \global\ladderCount=0 \let\next= \dLRTotherArrows\def\xLRTRungsx{#5\\} \else%
\let\next= \dLRTRungs\def\xLRTRungsx{#4\\#5\\}\fi\expandafter\next\xLRTRungsx%
}

\def \dLRTotherArrows#1#2#3\\{%
\global\advance\ladderCount by 1 %
\nextRungLRT\ladderCount\advance\nextRungLRT by 1 %
\ifnum\ladderType=1 %
\def\xdLRTotherArrowsx{#1}%
\ifx\xdLRTotherArrowsx\LRTarrowKill \else\DhA[\LRToptionName]{\\{\number\ladderCount 1}{\number\nextRungLRT 1}{#1}}\fi%
\def\xdLRTotherArrowsx{#2}%
\ifx\xdLRTotherArrowsx\LRTarrowKill \else\DhA[\LRToptionName]{\\{\number\ladderCount 2}{\number\nextRungLRT 2}{#2}}\fi%
\else
\def\xdLRTotherArrowsx{#1}%
\ifx\xdLRTotherArrowsx\LRTarrowKill \else\DvA[\LRToptionName]{\\{1\number\ladderCount}{1\number\nextRungLRT}{#1}}\fi%
\def\xdLRTotherArrowsx{#2}%
\ifx\xdLRTotherArrowsx\LRTarrowKill \else\DvA[\LRToptionName]{\\{2\number\ladderCount}{2\number\nextRungLRT}{#2}}\fi%
\fi
\ifnum\nextRungLRT=\rungs \let\next=\relax %
\def\xdLRTotherArrowsx{}\else\let\next= \dLRTotherArrows\def\xdLRTotherArrowsx{#3\\}\fi%
\expandafter\next\xdLRTotherArrowsx
}

\def\LRTsetLRToption#1#2#3@{\gdef\LRToption{[#1]{#2}{-#3}}\gdef\LRToptionName{#1}%
\expandafter\setCoor\LRToption%
}%

\newenvironment{hLadder}[2][{}{0}{0}]{%
\global\ladderType=1 %
\global\ladderCount=1 %
\gdef\deltaXlist{\\1{0}}%
\gdef\deltaYlist{\\1{0}\\2{-#2}}%
\expandafter\LRTsetLRToption#1@%
\ifnum\tikzLRTdepthCount=0 %
\tikizSetup%
\fi\global\advance\tikzLRTdepthCount by 1\dLRTX
}
{%
\global\advance\tikzLRTdepthCount by -1 \ifnum\tikzLRTdepthCount=0 %
\gdef\tikizSetup{\begin{tikzpicture}}%
\gdef\arrowData{[->]}%
\end{tikzpicture}\fi%
}

\newenvironment{vLadder}[2][{}{0}{0}]{%
\global\ladderType=0 %
\global\ladderCount=1 %
\gdef\deltaYlist{\\1{0}}%
\gdef\deltaXlist{\\1{0}\\2{#2}}%
\expandafter\LRTsetLRToption#1@%
\ifnum\tikzLRTdepthCount=0 %
\tikizSetup%
\fi\global\advance\tikzLRTdepthCount by 1\dLRTX
}
{%
\global\advance\tikzLRTdepthCount by -1 \ifnum\tikzLRTdepthCount=0 %
\gdef\tikizSetup{\begin{tikzpicture}}%
\gdef\arrowData{[->]}%
\end{tikzpicture}\fi%
}

\renewcommand{\raisedLabel}[1]{}

\usepackage{amsrefs}

\begin{document}

\def\newCommand#1#2{\expandafter\let\csname #1\endcsname=#2}
\def\cs#1{\csname #1\endcsname}
\def\loop{\Omega}
\def\bp#1{\ast_{_{#1}}}
\def\pe{p}
\def\linc#1{\mathfrak i_{#1}}
\def\WW{w}

\def\flp{\epsilon}

\def\loop{\Omega}
\def\pathslrt#1{P#1}
\newCommand{based path space}{\lrtpaths}
\newCommand{classifying map}{\WW}
\newCommand{base point}{\bp}
\newCommand{fibration projection}{\pe}
\newCommand{loop inclusion}{\linc}
\def\lrtfl#1{L\hskip1pt#1}
\newCommand{free loops}{\lrtfl}
\newCommand{free loop projection}{\flp}

\newCommand{bi-torsor set}{T}
\newCommand{bi-torsor left group}{G}
\newCommand{bi-torsor right group}{H}
\newcommand{\lrtrtlI}[1][]{\bar{\gamma}_{#1}}
\newCommand{right to left isomorphism}{\lrtrtlI}
\newCommand{left to right isomorphism}{\gamma}

\newcommand{\PA}[2]{#1\bullet #2}
\newCommand{path addition}{\PA}
\newcommand{\bps}[1]{P\hskip1pt#1}
\newCommand{based path space}{\bps}

\newcommand{\bls}[1]{\Omega\hskip1pt #1}
\newCommand{based loop space}{\bls}

\newCommand{principal fibration projection}{p}

\def\cmd{B}
\newCommand{classifying map domain}{\cmd}
\def\cmr{C}
\newCommand{classifying map range}{\cmr}

\def\cs#1{\expandafter\ifx\csname #1\endcsname\relax
\errmessage{no command #1}%
\def\xx{\hbox{\textcolor{red}{\bf??}}}\else%
\def\xx{\csname #1\endcsname}\fi\xx}

\def\LLifts#1#2#3#4#5{%
Lift^{\scriptscriptstyle#2}_{\ast}\bigl(\hskip-3pt\lower3pt\hbox{$\begin{tikzpicture}
\node (A) at (0,0) {$#1$};
\node (B) at (1.2,0) {$#5$};
\draw [->,dotted] (A) to node [above] {} (B);
\end{tikzpicture}$}\hskip-3pt\bigr) }
\newCommand{Lifts}{\LLifts}

\newCommand{free loops section}{s}

\def\lrtlh#1{{\mathcal L}#1_{\#}}
\newCommand{lift homomorphism}{\lrtlh}

\newcommand{\cma}[2]{\mathfrak c_{_{#1->#2}}}
\newCommand{constant map at}{\cma}
\newcommand{\fsA}[2]{{#1}^{#2}}
\newCommand{function space}{\fsA}
\newcommand{\bfsA}[2]{{#1}_{\ast}^{#2}}
\newCommand{based function space}{\bfsA}

\def\lrtcc{\mathcal K}
\newCommand{convenient category}{\lrtcc}

\newcommand{\piooo}[1][]{{#1}^{-1}}
\newCommand{path inverse}{\piooo}

\newcommand{\mfc}{\mathfrak c}
\newCommand{base point in function space}{\mfc}

\newcommand{\imofsX}[1]{\cs{classifying map}^{#1}}
\newCommand{induced map on function spaces}{\imofsX}

\def\lrtsqm{\mathbf s}
\newCommand{smash quotient map}{\lrtsqm}

\newcommand{\lrtls}[1][]{\Omega^{#1}\Sigma}
\newCommand{loop suspension}{\lrtls}

\def\coker#1{{\rm coker}\left(#1\right)}
\def\lrtcc#1{P_{#1}}
\newCommand{cohomology cokernel}{\lrtcc}

\def\lrtpm#1{\mathfrak s_{#1}}
\newCommand{power map}{\lrtpm}
\def\lrtpl#1{\psi^{\prime}_{#1}}
\newCommand{power lift}{\lrtpl}

\def\lrtcs#1{{\mathbf B}#1}
\newCommand{classifying space}{\lrtcs}
\def\lrtcps#1{\mathbb{CP}^{#1}}
\newCommand{CP}{\lrtcps}
\def\lrthps#1{\mathbb{HP}^{#1}}
\newCommand{HP}{\lrthps}

\newcommand{\lrthtxs}[1][]{\psi_{#1}}
\newCommand{homomorphism to [X,S^3]}{\lrthtxs}

\newCommand{classifying map total space}{E}

\title{The principal fibration sequence and the second cohomotopy set}
\author{Laurence R. Taylor}
\address{Department of Mathematics\\University of Notre Dame\\Notre Dame, IN 46556}
\email{taylor.2@nd.edu}
\begin{abstract}
Let $\cs{principal fibration projection}\colon \cs{classifying map total space} 
\to \cs{classifying map domain}$ 
be a principal fibration with classifying map 
$\cs{classifying map}\colon \cs{classifying map domain} \to \cs{classifying map range}$.
It is well-known that the group
$[X,\cs{based loop space}{\cs{classifying map range}}]$ acts on 
$[X,\cs{classifying map total space}]$ with orbit space the image of 
$\cs{principal fibration projection}_{\#}$, where 
$\cs{principal fibration projection}_{\#}\colon [X,\cs{classifying map total space}] 
\to [X,\cs{classifying map domain}]$. 
The isotropy subgroup of the map of $X$ to the base point of 
$\cs{classifying map total space}$ is 
also well-known to be the image of $[X,\cs{based loop space}\cs{classifying map domain}]$. 
The isotropy subgroups for other maps $e\colon X \to 
\cs{classifying map total space}$ 
can definitely change as $e$ does. 

The set of homotopy classes
of lifts of $f$ to the free loop space on $\cs{classifying map domain}$ is a group. 
If $f$ has a lift to $\cs{classifying map total space}$,  
the set $\cs{principal fibration projection}_{\#}^{-1}(f)$ 
is identified with the cokernel of a natural homomorphism 
from this group of lifts to $[X,\cs{based loop space}{\cs{classifying map range}}]$.  

As an example, $[X,S^2]$ is enumerated for $X$ a 4-complex. 
\end{abstract}
\maketitle

\section{Results and discussion}
For based spaces $X$ and $Y$, let $[X,Y]$ denote the set of based homotopy classes 
of maps from $X$ to $Y$. 
The constant map to the base point makes $[X,Y]$ into a based set. 
If $Y$ is based, the constant path at the base point is a base point for the
based loops, $\cs{based loop space}Y$.

A principal fibration 
$\cs{principal fibration projection}\colon \cs{classifying map total space} \to \cs{classifying map domain}$ 
is a fibration with a classifying map
$\cs{classifying map}\colon \cs{classifying map domain} \to \cs{classifying map range}$ 
such that $\cs{classifying map total space}$ is a pull-back of the path-loop
fibration for $\cs{classifying map range}$ along $\cs{classifying map}$. 
Pick a base point $\cs{base point}{\cs{classifying map total space}} \in 
\cs{classifying map total space}$. 
Let the base point in $\cs{classifying map domain}$ be
$\cs{base point}{\cs{classifying map domain}} = 
\cs{principal fibration projection}(\cs{base point}{\cs{classifying map total space}})$ and
let the base point in $\cs{classifying map range}$ be 
$\cs{base point}{\cs{classifying map range}} = 
\cs{classifying map}(\cs{base point}{\cs{classifying map domain}})$, so that 
$\cs{classifying map}$ and $\cs{principal fibration projection}$ 
become based maps. 

It is a result going back to Peterson \cite{Peterson}*{Lemma 2.1, p.~199} and 
Nomura \cite{Nomura}*{Corollary  2.1, p.~118}
that there is an exact sequence 
of based sets
\newNumber{fibre sequence}
\[
\cdots \to [X, \cs{based loop space}\cs{classifying map domain}] 
\to [X, \cs{based loop space}{\cs{classifying map range}}] 
\to [X, \cs{classifying map total space}] 
\rightlabeledarrow{\cs{principal fibration projection}_{\#}}{} 
[X,\cs{classifying map domain}] \rightlabeledarrow{\cs{classifying map}_{\#}}{} 
[X,\cs{classifying map range}]
\leqno(\ref{fibre sequence})\]
in the following sense. 
Each map is a map of based sets 
and the image of one map is
the inverse image of the base point for the following map.

One way to derive this sequence is to fix an 
$f\colon X \to \cs{classifying map domain}$ and 
consider the set 
$\cs{Lifts}X f {\cs{classifying map domain}} {\cs{principal fibration projection}}
{\cs{classifying map total space}}$ 
of based homotopy classes of lifts of $f$. 
There is a forgetful map 
$\cs{Lifts}X f {\cs{classifying map domain}} {\cs{principal fibration projection}}
{\cs{classifying map total space}}
\to [X,\cs{classifying map total space}]$ and properties of
fibrations imply that the image is 
$\cs{principal fibration projection}_{\#}^{-1}(f)$. 

Peterson and Thomas \cite{Peterson-Thomas}*{Lemma 4.1, p.~17} show 
that there is a left action of the group 
$[X, \cs{based loop space}{\cs{classifying map range}}]$ on the set 
$[X,\cs{classifying map total space}]$ 
which identifies the orbit space with 
$\cs{classifying map}_{\#}^{-1}(\ast)$ where $\ast\in [X,\cs{classifying map range}]$ 
is the base point.
They show that the set 
$\cs{Lifts}X f {\cs{classifying map domain}} {\cs{principal fibration projection}}
{\cs{classifying map total space}}$ 
is a left $[X,\cs{based loop space}{\cs{classifying map range}}]$ torsor and 
this gives exactness in (\ref{fibre sequence}) at 
$[X,\cs{based loop space}{\cs{classifying map range}}]$. 
It further follows that the isotropy subgroup of the action on the base point 
of $[X,\cs{classifying map total space}]$ 
is the image of $[X, \cs{based loop space}{\cs{classifying map domain}}]$. 

\bigskip

There is another way to proceed.
For a space $Y$, let $\cs{free loops}Y$ denote the free loop space and
let $\cs{free loop projection}\colon \cs{free loops}Y \to Y$ denote the
projection.
The map $\cs{free loop projection}$ is a fibration. 
The constant loop at $y\in Y$ defines a section
$\cs{free loops section}\colon Y \to \cs{free loops}Y$ 
so the set of lifts has a base point, $s\circ f$. 
If $\cs{base point}{Y}\in Y$ is a base point, the space $\cs{free loops}Y$ has a base point, 
$\cs{free loops section}(\cs{base point}{Y})$. 
Let  $\cs{Lifts}X f Y {\cs{free loop projection}} {\cs{free loops}Y}$ denote
the set of based homotopy classes of lifts of $f$.

Addition of loops makes $\cs{Lifts}X f Y {\cs{free loop projection}} {\cs{free loops}Y}$
into a group. 
A based map 
$\cs{classifying map}\colon \cs{classifying map domain}\to \cs{classifying map range}$ 
induces a based map 
$\cs{free loops}{\cs{classifying map}} 
\colon \cs{free loops}\cs{classifying map domain} \to 
\cs{free loops}\cs{classifying map range}$ and a family of 
group homomorphisms 
\[
\cs{Lifts}X f {\cs{classifying map domain}} 
{\cs{free loop projection}}{\cs{free loops}\cs{classifying map domain}}
\rightlabeledarrow{\ \cs{lift homomorphism}{\cs{classifying map}}\ }{}
\cs{Lifts}X {g \circ f} {\cs{classifying map range}} 
{\cs{free loop projection}}{\cs{free loops}\cs{classifying map range}}
\]

\begin{Theorem}\namedLabel{isotropy subgroup}
The set $\cs{Lifts}X f {\cs{classifying map domain}} {\cs{principal fibration projection}}
{\cs{classifying map total space}}$ is a right 
$\cs{Lifts}X{\cs{classifying map}\circ f} {\cs{classifying map range}}
{\cs{free loop projection}} {\cs{free loops}{\cs{classifying map range}}}$ torsor. 
To each element $e \in [X,\cs{classifying map total space}]$ 
lifting $f \in [X, \cs{classifying map domain}]$, there is associated 
a group isomorphism 
\[
\cs{Lifts}X{\cs{classifying map}\circ f} {\cs{classifying map range}} 
{\cs{free loop projection}} {\cs{free loops}{\cs{classifying map range}}} 
\rightlabeledarrow{\ \cs{right to left isomorphism}[e]\ }{}  
[X,\cs{based loop space}{\cs{classifying map range}}]
\]
The image of the composition
\[
\cs{Lifts}X {f} {\cs{classifying map domain}} 
{\cs{free loop projection}} {\cs{free loops}{\cs{classifying map domain}}} 
\rightlabeledarrow{\ \cs{lift homomorphism}{\cs{classifying map}}\ }{} 
\cs{Lifts}X{\cs{classifying map}\circ f} {\cs{classifying map range}} 
{\cs{free loop projection}} {\cs{free loops}{\cs{classifying map range}}} 
\rightlabeledarrow{\ \cs{right to left isomorphism}[e]\ }{}  
[X,\cs{based loop space}{\cs{classifying map range}}]
\]
is the isotropy subgroup of $e$ under the left 
$[X,\cs{based loop space}{\cs{classifying map range}}]$ action on 
$[X, \cs{classifying map total space}]$. 
\end{Theorem}

\namedRef{isotropy subgroup} gives a sequence with 
many of the same properties as  (\ref{fibre sequence}). 
The group $[X,\cs{based loop space}{\cs{classifying map range}}]$ 
is the same for both sequences. 
The group $\cs{Lifts}X {f} {\cs{classifying map domain}} 
{\cs{free loop projection}} {\cs{free loops}{\cs{classifying map domain}}}$ 
depends on $\cs{classifying map domain}$, $f$ and 
$X$ instead of just $X$ and $\cs{based loop space}{\cs{classifying map domain}}$ 
as in (\ref{fibre sequence}), but it is still independent of ${\cs{classifying map range}}$ and $\cs{classifying map}$. 
The homomorphism between these two groups can depend on $f$ 
in addition to just $\cs{classifying map}$ (see \S\ref{second cohomotopy set}). 
The additional information supplied by \namedRef{isotropy subgroup} 
comes from the fact that all the isotropy subgroups are determined 
rather than just the isotropy subgroup 
of the null homotopic map as in (\ref{fibre sequence}). 

\begin{Remark}
J.~Rutter \cite{Rutter} has results similar to these if $\cs{classifying map domain}$ 
and $\cs{classifying map range}$ are H-spaces. 
In this case the multiplication can be used to naturally identify
$\cs{Lifts}X {f} {\cs{classifying map domain}} 
{\cs{free loop projection}} {\cs{free loops}{\cs{classifying map domain}}}$ with
$[X,\cs{based loop space}{\cs{classifying map domain}}]$. 
Rutter uses the H-space structure to describe a homomorphism
$[X,\cs{based loop space}{\cs{classifying map domain}}] \to
[X,\cs{based loop space}{\cs{classifying map range}}]$, depending on $f$, 
which presumably is related to the homomorphism given by
\namedRef{isotropy subgroup} whenever $B$ is an H-space. 
In general this homomorphism can not be the one induced by 
$\cs{based loop space}{\cs{classifying map}}$ since the size of the cokernel 
can vary with $f$. (See \ref{Pontrjagin's calculation} and 
\ref{4-dimensional example}.) 

The observation that there is a natural right action of the section groups 
adds generality, and perhaps clarity, to the result. 
\end{Remark}

\medskip

An additional observation is that the calculations required by 
(\ref{fibre sequence}) and \namedRef{isotropy subgroup} are
natural in both the space $X$ and the principal fibration.

\section{Recall of some basic results}

The sequence (\ref{fibre sequence}) 
can be derived from standard results
about the path-groupoid applied to function spaces. 
The needed results are recalled below.
To prove \namedRef{isotropy subgroup} requires an additional 
technical lemma,
\namedRef{half free homotopy}. 

\subsection{Point set topology}
As usual all constructions take place in a ``convenient category'', 
$\cs{convenient category}$.
Vogt \cite{V} is a good reference.
One key point is that the \emph{exponential correspondence} holds,
the space of maps $X\times Y$ to $W$, is homeomorphic to 
the space of maps of $X$ to $\cs{function space}{W}{Y}$. 
Here the product gets the product topology in $\cs{convenient category}$
and $\cs{function space}{W}{Y}$ gets the topology given by 
starting with the compact-open topology and making it compactly-generated. 
Also, the \emph{subspace topology} on a subset is the one given by taking the
usual subspace topology and then making it compactly-generated. 

If $W_0\subset W$ is a subspace, in the category $\mathcal K$, 
then $\cs{function space}{W_0}{Y}$ with its topology is a subspace
of $\cs{function space}{W}{Y}$ with its topology. 

Given any point $w\in W$ and any space $Y$, 
let $\cs{constant map at}{Y}{w} \in \cs{function space}{W}{Y}$ 
denote the constant map of $Y$ to $w$. 
Anytime $W$ has a base point $\cs{base point}W \in W$, 
the map $\cs{constant map at}{Y}{\cs{base point}W}$ will be 
the base point in $\cs{function space}{W}{Y}$. 
If both $Y$ and $W$ are based, then
$\cs{based function space}{W}{Y}$ is the subspace of 
$\cs{function space}{W}{Y}$ consisting of all maps $f\colon Y \to W$
which preserve the base points. 

A base point is non-degenerate provided the pair 
$(W,\cs{base point}W)$ is an NDR pair. 

If $(W,W_0)$ and $(Y,Y_0)$ are pairs, $\cs{function space}{(W,W_0)}{(Y,Y_0)}$ 
denotes the space of all continuous functions $Y \to W$ sending $Y_0 \to W_0$.
It is given the subspace topology in $\mathcal K$ from $\cs{function space}{W}{Y}$.

\begin{Result}
If $(W,W_0)$ is an NDR pair and if $Y$ is compact, then 
$\cs{function space}{W_0}{Y}$ is a subspace of 
$\cs{function space}{(W,W_0)}{(Y,Y_0)}$ and the pair is an NDR pair. 
\end{Result}
\begin{proof}
Since a subspace of a subspace is a subspace 
$\cs{function space}{W_0}{Y}$ is a subspace of 
$\cs{function space}{(W,W_0)}{(Y,Y_0)}$. 

If $u\colon W \to [0,1]$ is the map which is part of the definition of an NDR pair,
then $\hat{u}\colon \cs{function space}{W}{Y} \to [0,1]$ defined by
$\hat{u}(f) = {\rm sup}_{y\in Y}\ u\bigl(f(y)\bigr)$ is continuous. 
This uses $Y$ compact. 
Note $\hat{u}^{-1}(0) = \cs{function space}{W_0}{Y}$. 

If $F\colon W\times [0,1] \to W$ is the homotopy which is the other 
part of the definition of an NDR pair, then if
$\hat{F}\colon \cs{function space}{W}{Y} \times [0,1] \to \cs{function space}{W}{Y}$ 
is defined by $\bigl(\hat{F}(f,t)\bigr)(y) = F\bigl(f(y),t\bigr)$, the pair 
$(\hat{F},\hat{u})$ shows the function spaces form an NDR pair.  
\end{proof}
\begin{Remark}
If $\cs{base point}W\in W$ is non-degenerate then 
$\cs{constant map at}{Y}{\cs{base point}W}$ is a non-degenerate point in both 
$\cs{function space}{W}{Y}$ and $\cs{based function space}{W}{Y}$. 
\end{Remark}

\subsection{The path groupoid}
Given two points $w_0$, $w_1\in W$ let 
$W_{w_0,w_1}$ denote the set of homotopy classes of paths 
from $w_0$ to $w_1$ where the homotopies are rel end points.
The set $W_{w_0,w_1}$ is non-empty if and only if $w_0$ and
$w_1$ are in the same path component of $W$. 

If $w_0$, $w_1$ and $w_2$ are all in one path component of $W$,
path concatenation defines an associative pairing
\[
W_{w_0,w_1} \times W_{w_1,w_2} \to W_{w_0,w_2}
\]
Reversing the path defines an involution  
$W_{w_0,w_1} \rightlabeledarrow{\ \cs{path inverse}\ }{}
W_{w_1,w_0}$, and hence a bijection, 
such that the image of the composition 
$W_{w_0,w_1} \rightlabeledarrow{\ 1\times \cs{path inverse}\ }{}
W_{w_0,w_1} \times W_{w_1,w_0} 
\to
W_{w_0,w_0}$ is the constant path at $w_0$. 
There is a similar constant map $W_{w_0,w_1} \to W_{w_1,w_1}$. 

For any $w\in W$, $W_{w,w}$ is a group under path concatenation with 
$\cs{path inverse}$ being the inverse map. 

If $W_{w_0,w_1}$ is non-empty, the group $W_{w_0,w_0}$ acts
on it on the left and the group $W_{w_1,w_1}$ acts on the right.
Associativity of path concatenation makes 
$W_{w_0,w_1}$ into a bi-set. 

\begin{Result}\namedLabel{general torsor}
If non-empty, the bi-set $W_{w_0,w_1}$ is a torsor for each group action.
\end{Result}
\begin{proof}
To be a torsor means the group action is transitive and the isotropy 
subgroup of any point is trivial. 

Let $\tau_0$, $\tau_1 \in W_{w_0,w_1}$.
Then $\tau_1 = \tau_0\cs{path addition}{}(\tau_0^{-1}\cs{path addition}{}\tau_1)$ 
and $\tau_0^{-1}\cs{path addition}{}\tau_1 \in W_{w_1,w_1}$.
Similarly
$\tau_1 = (\tau_1\cs{path addition}{} \tau_0^{-1})\cs{path addition}{}\tau_0$ 
and $\tau_1\cs{path addition}{} \tau_0^{-1} \in W_{w_0,w_0}$. 
Hence both actions are transitive. 

Now suppose $\lambda \cs{path addition}{}\tau = \tau$ 
for some $\lambda \in W_{w_0,w_0}$ and some $\tau \in W_{w_0,w_1}$.
Then $(\lambda \cs{path addition}{}\tau)\cs{path addition}{} \tau^{-1} = 
\tau\cs{path addition}{} \tau^{-1}$ and therefore
$\lambda$ is homotopic rel end points to the constant path and so 
the isotropy subgroup of $\tau$ under the left action is trivial. 
A similar calculation shows the right action also has trivial isotropy subgroups.
\end{proof}

\begin{Lemma}\namedLabel{half free homotopy}
Let $\tau_0$, $\tau_1$ be representatives of elements in $W_{w_0,w_1}$ and 
let $\phi\in W_{w_1,w_1}$. 
There exists a homotopy
\[F\colon [0,1]\times [0,1] \to W
\]
with $F(t,0) = \tau_0(t)$, $F(t,1) = \tau_1(t)$, 
$F(0,s) = w_0 = \tau_0(0)=\tau_1(0)$ and $F(1,s) = \phi(s)$
if and only if 
$\tau_1 = \tau_0 \cs{path addition}{} \phi$.
\end{Lemma}
\begin{proof}
Let $H$ be a homotopy rel end points from $\tau_0 \cs{path addition}{} \phi$ 
to $\tau_1$. 
Figure 1 is a visual representation for $H$ and Figure 2 is one for $F$. 
\namedRef{half free homotopy} is equivalent to constructing $F$ given $H$ and $H$ given $F$.

\def\LLx{0}\def\LLy{0}\def\LL{1.75}\def\HH{1.75}

\begin{tikzpicture}[inner sep=.5mm,]
\node (ll) at (\LLx,\LLy) [circle,draw, fill=black] {};
\node (lr) at (\LLx+\LL,\LLy) [circle,draw, fill=black] {};
\node (ul) at (\LLx,\LLy+\HH) [circle,draw, fill=black] {};
\node (half) at (\LLx,\LLy+\HH/2) [circle,draw, fill=black] {};
\node (ur) at (\LLx+\LL,\LLy+\HH) [circle,draw, fill=black] {};
\node (W) at (\LLx+\LL+2.6,\LLy+\HH/2) {$W$};
\draw[-](ll) to node[above=2pt]{$\scriptstyle w_0$} (lr);
\draw[-](ll) to node[left=1pt]{$\tau_0$} (half);
\draw[-](half) to node[left=6pt]{$\phi$} (ul);
\draw[-](ul) to node[above=2pt]{$\scriptstyle w_1$} (ur);
\draw[-] (lr) to node[left=1pt]{$\tau_1$} (ur);
\draw[->] (\LLx+\LL+.2,\LLy+\HH/2) parabola [parabola height=4mm] +(2,0);
\path (\LLx+\LL+.2,\LLy+\HH/2) to node [above=14pt] {$H$} (W);
\node at (\LLx +\LL/2,\LLy-20pt) {Figure 1};
\end{tikzpicture}
\hskip .5in
\begin{tikzpicture}[inner sep=.5mm,]
\node (ll) at (\LLx,\LLy) [circle,draw, fill=black] {};
\node (lr) at (\LLx+\LL,\LLy) [circle,draw, fill=black] {};
\node (ul) at (\LLx,\LLy+\HH) [circle,draw, fill=black] {};
\node (ur) at (\LLx+\LL,\LLy+\HH) [circle,draw, fill=black] {};
\node (W) at (\LLx+\LL+2.6,\LLy+\HH/2) {$W$};
\draw[-](ll) to node[above=2pt]{$\scriptstyle w_0$} (lr);
\draw[-](ll) to node[left=1pt]{$\tau_0$} (ul);
\draw[-](ul) to node[above=2pt]{$\phi$} (ur);
\draw[-] (lr) to node[left=1pt]{$\tau_1$} (ur);
\draw[->] (\LLx+\LL+.2,\LLy+\HH/2) parabola [parabola height=4mm] +(2,0);
\path (\LLx+\LL+.2,\LLy+\HH/2) to node [above=14pt] {$F$} (W);
\node at (\LLx +\LL/2,\LLy-20pt) {Figure 2};
\end{tikzpicture}

\hskip 1.25in 
\begin{tikzpicture}[inner sep=.5mm,]
\node (ll) at (\LLx,\LLy) [circle,draw, fill=black] {};
\node (lr) at (\LLx+\LL,\LLy) [circle,draw, fill=black] {};
\node (ur) at (\LLx+\LL,\LLy+\HH) [circle,draw, fill=black] {};
\node (half) at (\LLx + \LL/2,\LLy + \HH/2) [circle,draw, fill=black] {};
\node (W) at (\LLx+\LL+2.6,\LLy+\HH/2) {$W$};
\draw[-](ll) to node[above=2pt]{$\scriptstyle w_0$} (lr);
\draw[-] (lr) to node[left=0pt]{\lower 16pt\hbox{$\tau_1$}} (ur);
\draw[-] (ll) tonode[left=2pt]{$\tau_0$} (half);
\draw[-] (half) tonode[left=2pt]{$\phi$} (ur);
\draw[->] (\LLx+\LL+.2,\LLy+\HH/2) parabola [parabola height=4mm] +(2,0);
\path (\LLx+\LL+.2,\LLy+\HH/2) to node [above=14pt] {$G$} (W);
\node at (\LLx +\LL/2,\LLy-20pt) {Figure 3};
\end{tikzpicture}

There is an evident map from the squares in Figures 1 and 2 to 
the triangle in Figure 3. 
Either map $H$ or $F$ induces a map $G$ from the triangle to $W$. 
Given $G$, composition with the map from the appropriate square constructs 
both $F$ and $H$.
\end{proof}

\subsection{Bi-torsors}
Suppose $\cs{bi-torsor set}$ is a left $\cs{bi-torsor left group}$ torsor 
and a right $\cs{bi-torsor right group}$ torsor as well as 
a $\cs{bi-torsor left group}$-$\cs{bi-torsor right group}$ biset.
For $x\in \cs{bi-torsor set}$, 
define functions $\cs{left to right isomorphism}_x\colon 
\cs{bi-torsor left group} \to \cs{bi-torsor right group}$ and 
$\cs{right to left isomorphism}[x]\colon 
\cs{bi-torsor right group} \to \cs{bi-torsor left group}$
by $g\bullet x = x\bullet \cs{left to right isomorphism}_x(g)$ and 
$\cs{right to left isomorphism}[x](h)\bullet x = x \bullet h$.

\begin{Result}\namedLabel{biset isomorphisms}
Each $x\in \cs{bi-torsor set}$ defines a group isomorphism
$\cs{left to right isomorphism}_x\colon
\cs{bi-torsor left group} \to \cs{bi-torsor right group}$ 
and an inverse isomorphism 
$\cs{right to left isomorphism}[x]\colon
\cs{bi-torsor right group} \to \cs{bi-torsor left group}$
\end{Result}
\begin{proof}
Note $\cs{left to right isomorphism}_x
(e_{\cs{bi-torsor left group}}) = e_{\cs{bi-torsor right group}}$ and
$\cs{right to left isomorphism}[x]
(e_{\cs{bi-torsor right group}}) = e_{\cs{bi-torsor left group}}$. 
Also check 
$x\bullet \cs{left to right isomorphism}_x(g_1 g_2) = (g_1 g_2)\bullet x 
= g_1\bullet (g_2\bullet x) = g_1\bullet\bigl(x \bullet 
\cs{left to right isomorphism}_x(g_2)\bigr) =
( g_1\bullet x) \bullet \cs{left to right isomorphism}_x(g_2) =
x \bullet\bigl(\cs{left to right isomorphism}_x(g_1)
\cs{left to right isomorphism}_x(g_2)\bigr)
$
so $\cs{left to right isomorphism}_x$ is multiplicative. 
Check $\cs{left to right isomorphism}_x$ and 
$\cs{right to left isomorphism}[x]$ are inverse functions. 
Hence they are inverse homomorphisms and the result follows.
\end{proof}
\begin{Result}
If $x_1$, $x_2 \in \cs{bi-torsor set}$, then 
$\cs{left to right isomorphism}_{x_1}$
and $\cs{left to right isomorphism}_{x_2}$ are conjugate as are
$\cs{right to left isomorphism}[x_1]$
and $\cs{right to left isomorphism}[x_2]$
\end{Result}
\begin{proof}
If $x_1 = x_2 \bullet h$, 
$x_1 \bullet \cs{left to right isomorphism}_{x_1}(g) = g\bullet x_1 = 
g\bullet(x_2\bullet h) = 
x_2\bullet(\cs{left to right isomorphism}_{x_2}(g) h) = 
(x_2\bullet h)\bullet(h^{-1}
\cs{left to right isomorphism}_{x_2}(g) h) = 
x_1 \bullet(h^{-1}\cs{left to right isomorphism}_{x_2}(g) h)$ so
$\cs{left to right isomorphism}_{x_1}(g) = 
h^{-1} \cs{left to right isomorphism}_{x_2} (g) h$. 
The proof for the $\cs{right to left isomorphism}$ is similar. 
\end{proof}

\subsection{Principal Fibrations}
A \emph{principal fibration} is a fibration 
$\cs{classifying map total space}\rightlabeledarrow{p}{}  \cs{classifying map domain}$ 
which is a pull-back of the path-loop fibration 
$\Omega {\cs{classifying map range}} \to 
\cs{based path space}{\cs{classifying map range}} 
\rightlabeledarrow{\rho_{\cs{classifying map range}}}{} {\cs{classifying map range}}$
along a map 
$\cs{classifying map}\colon \cs{classifying map domain} \to {\cs{classifying map range}}$. 
The definition of the space $\cs{based path space}{\cs{classifying map range}}$ 
requires a base point in ${\cs{classifying map range}}$, 
say $\cs{base point}{{\cs{classifying map range}}}$.
Then $\cs{based path space}{\cs{classifying map range}}$ is the space of all maps 
$\lambda\colon [0,1] \to {\cs{classifying map range}}$ 
such that $\lambda(0)=\cs{base point}{{\cs{classifying map range}}}$. 
Equivalently it is the subspace of 
$\cs{function space}{{\cs{classifying map range}}}{[0,1]}$ of paths 
that start at $\cs{base point}{\cs{classifying map range}}$, 
$\cs{function space}{({\cs{classifying map range}},
\cs{base point}{{\cs{classifying map range}}})}{([0,1],0)}$.

Up to fibre homotopy equivalence, a principal fibration has a standard model.
The total space is 
$\cs{classifying map total space}_{\cs{classifying map}} \subset 
\cs{classifying map domain} \times 
\cs{function space}{\cs{classifying map range}}{[0,1]}$ 
such that $(b,\lambda) \in \cs{classifying map total space}_{\cs{classifying map}}$ 
if and only if $\cs{classifying map}(b) = \lambda(1)$
and $\cs{base point}{\cs{classifying map range}} = \lambda(0)$. 
The fibration projection is just projection onto the $\cs{classifying map domain}$ factor. 
If $\cs{classifying map domain}$ is given a base point 
$\cs{base point}{\cs{classifying map domain}}$ such that 
$\cs{classifying map}(\cs{base point}{\cs{classifying map domain}}) = 
\cs{base point}{\cs{classifying map range}}$, then
$\cs{classifying map total space}_{\cs{classifying map}}$ has a base point,
$\bigl(\cs{base point}{\cs{classifying map domain}}, \cs{base point in function space}\bigr)$ 
where it should cause no confusion to shorten the notation for the base 
point in a function space to $\cs{base point in function space}$. 

For the purposes of this paper it suffices to pick a convenient based map 
$\cs{classifying map}$, and then work with 
$\cs{classifying map total space}_{\cs{classifying map}}$. 
Two $\cs{classifying map}$ which are based homotopic yield 
$\cs{classifying map total space}_{\cs{classifying map}}$ which are based fibre homotopy
equivalent and all questions discussed here only depend on the based fibre 
homotopy type of the fibration. 

The next result describes the set of lifts. 
There is a map induced by composition with $\cs{classifying map}$,
$\cs{induced map on function spaces}{X}\colon \cs{function space}
{\cs{classifying map domain}}{X} \to 
\cs{function space}{{\cs{classifying map range}}}{X}$.

\begin{Result}
The set of homotopy classes of lifts of 
$f \in \cs{function space}{\cs{classifying map domain}}{X}$ 
is equivalent to the set
$W_{w_0,w_1}$ where $W = \cs{function space}{{\cs{classifying map range}}}{X}$, 
$w_0 = \cs{base point in function space}$ and 
$w_1 =  \cs{induced map on function spaces}{X}(f)$. 
If $f$ is based, then the set of based homotopy classes of lifts of 
$f$ is equivalent to the set
$W_{w_0,w_1}$ for the same $w_i$ but with 
$W = \cs{based function space}{{\cs{classifying map range}}}{X}$.
\end{Result}
\begin{proof}
A map of $X$ to 
$\cs{classifying map total space}_{\cs{classifying map}}$ consists of a map 
$f\colon X \to \cs{classifying map domain}$ 
and a map $\Lambda\colon X \to \cs{function space}{{\cs{classifying map range}}}{[0,1]}$ 
satisfying two conditions:
$\cs{classifying map}\circ f(x) = \Lambda(x,1)$ and 
$\Lambda(x,0) = \cs{base point}{\cs{classifying map range}}$. 

Consider the map $f$ as a point $f\in \cs{function space}
{\cs{classifying map domain}}{X}$ and the map 
$\Lambda$ as a map $\Lambda\colon [0,1] \to 
\cs{function space}{{\cs{classifying map range}}}{X}$ 
satisfying two conditions:
$\Lambda(0) = \cs{base point in function space}$
and $\Lambda(1) = \cs{induced map on function spaces}{X}(f)$. 
Two lifts of $f$, $\Lambda_0$ and $\Lambda_1$, are homotopic 
as lifts if and only if $\Lambda_0$ and $\Lambda_1$ are 
homotopic rel end-points, that is, they represent the same element in 
$W_{w_0,w_1}$.
\end{proof}

\begin{Result}
Given a map $f\in \cs{function space}{Y}{X}$, 
a lift to the free loop space is a map 
$\Phi\colon [0,1] \to \cs{function space}{Y}{X}$ such that $\Phi(0) = \Phi(1) = f$. 
In other words, the set of homotopy classes of lifts of $f$ to the free 
loop space on $Y$ is equivalent to $W_{w,w}$ with 
$W = \cs{function space}{Y}{X}$ and $w = f$. 
If $f\in\cs{based function space}{Y}{X}$ then the based lifts are equivalent
to $W_{w,w}$ with the same $w$ and $W = \cs{based function space}{Y}{X}$. 
\end{Result}

\section{The proof of \namedRef{isotropy subgroup}}
Fix a principal fibration  $\cs{principal fibration projection} \colon 
\cs{classifying map total space}_{\cs{classifying map}} \to \cs{classifying map domain}$, 
$\cs{classifying map}\colon \cs{classifying map domain} \to {\cs{classifying map range}}$. 
Fix a base point in 
$\cs{classifying map domain}$ and use 
its image to base ${\cs{classifying map range}}$.
This gives a preferred base point in 
$\cs{classifying map total space}_{\cs{classifying map}}$. 
Also fix a based  space $X$. 

Let $W = \cs{based function space}{\cs{classifying map range}}{X}$. 
Since ${\cs{classifying map range}}$ 
must have a base point to define $\cs{classifying map total space}_{\cs{classifying map}}$,
let $\cs{base point in function space}$ be the constant map of $X$
to the base point of ${\cs{classifying map range}}$.
Fix $e\colon X \to \cs{classifying map total space}$ and let 
$f = \cs{principal fibration projection}\circ e$. 

\begin{Remark}
Given $f\colon X \to \cs{classifying map domain}$, 
there exist such $e$'s if and only if
$\cs{classifying map}\circ f$ is null-homotopic rel base point. 
\end{Remark}

Up to homotopy of lifts, $e$ is determined by $f$ and 
$\Lambda\in W_{\cs{base point in function space}, 
\cs{induced map on function spaces}{X}(f)}$. 
The group acting on the left is 
$W_{\cs{base point in function space},\cs{base point in function space}} = 
[X,\cs{based loop space}{\cs{classifying map range}} ]$. 
The group acting on the right is 
$W_{\cs{induced map on function spaces}{X}(f),
\cs{induced map on function spaces}{X}(f)} = 
\cs{Lifts}X {\cs{classifying map}\circ f} {\cs{classifying map domain}} 
{\cs{principal fibration projection}} {\cs{free loops}{\cs{classifying map range}}}$. 
The isomorphism 
$\cs{right to left isomorphism}[e]$ in 
\namedRef{isotropy subgroup} is the map defined by \namedRef{biset isomorphisms}. 

Two lifts $\Lambda_0$ and $\Lambda_1$ are homotopic in 
$[X,\cs{classifying map total space}]$ 
if and only if 
there are homotopies $\Phi\colon [0,1] \to 
\cs{based function space}{\cs{classifying map domain}}{X}$ 
with $\Phi(0)  = \Phi(1) = f$ and 
$F\colon [0,1]\times [0,1] \to 
\cs{based function space}{\cs{classifying map range}}{X}$ such that
$F(1,s) = \cs{induced map on function spaces}{X}\bigl(\Phi(s)\bigr)$, 
$F(i,t) = \Lambda_i(t)$, $i = 0$, $1$. 

Equivalently, 
$\Phi \in
\cs{Lifts}X {f} {\cs{classifying map domain}} 
{\cs{principal fibration projection}} {\cs{free loops}{\cs{classifying map domain}}}$
and, if 
$\phi = \cs{induced map on function spaces}{X}(\Phi) \in 
\cs{Lifts}X {\cs{classifying map}\circ f} {\cs{classifying map domain}} 
{\cs{principal fibration projection}} {\cs{free loops}{\cs{classifying map range}}}$, 
\namedRef{half free homotopy} completes the proof of 
\namedRef{isotropy subgroup}. 

\section{Some general remarks on calculations}
There are some situations in which the group of lifts calculation 
can be replaced by just calculating a set of homotopy classes of maps.

One situation, \namedRef{abelian lifts}, is a generalization of a result of James and Thomas, 
\cite{James-Thomas}*{Theorem 2.6, p.~493}. 
\begin{Theorem} 
Let $Y$ be a based space and let $f\colon X \to Y$ be a based map. 
Then
\[
\mathcal G = \cs{Lifts}{X}{f}{Y}{\cs{free loop projection}}{\cs{free loops}Y} 
\rightlabeledarrow{\ \iota\ }{}
[X, \cs{free loops}Y]
\rightlabeledarrow{\ \cs{free loop projection}_{\#}\ }{}
[X, Y]
\]
is exact in that the image of $\iota$ is
$\cs{free loop projection}_{\#}^{-1}(f)$. 
The image of $\iota$ is also the set of conjugacy classes
of elements of $\mathcal G$.
\end{Theorem}
\begin{proof}
A lift is a map $X \to \cs{free loops}Y = Y^{S^1}$. 
By the exponential correspondence a lift is also a map $S^1 \to Y^X$. 
The lift property is equivalent to the additional condition that 
the base point of $S^1$ goes to $f\in Y^X$. 
Hence $\mathcal G = \pi_1(Y^X; f)$. 

An element $[X, \cs{free loops}Y]$ is equal to an element in $[S^1, Y^X]$ 
with no condition on the base points except that the base point of $S^1$
lands in the path component of $f$. 
There is always a homotopy which takes the base point of $S^1$ to 
$f\in Y^X$ so the image of $\iota$ is
$\cs{free loop projection}_{\#}^{-1}(f)$. 

It is always true that the relation between $\pi_1(Y^X; f)$ and
the free homotopy classes is that the set of free homotopy classes 
is the set of conjugacy classes. 
\end{proof}
\begin{Corollary}[James \& Thomas, \cite{James-Thomas}]\namedLabel{abelian lifts}
The group 
$\cs{Lifts}{X}{f}{Y}{\cs{free loop projection}}{\cs{free loops}Y}$ is abelian 
if and only if $\iota$ is injective. 
\end{Corollary}

Given a map $\cs{classifying map}\colon \cs{classifying map domain} \to 
\cs{classifying map range}$, there is an induced map 
$\cs{free loops}\cs{classifying map}\colon \cs{free loops}\cs{classifying map domain} \to 
\cs{free loops}\cs{classifying map range}$
and
\[
\begin{matrix}
\cs{Lifts}{X}{f}{\cs{classifying map domain}}{\cs{free loop projection}}
{\cs{free loops}\cs{classifying map domain}} &
\rightlabeledarrow{\ \iota^1\ }{}&
[X, \cs{free loops}\cs{classifying map domain}]&
\rightlabeledarrow{\ \cs{free loop projection}^1_{\#}\ }{}&
[X, \cs{classifying map domain}]
\\
\downlabeledarrow[\big]{\cs{lift homomorphism}{\cs{classifying map}}}{}&&\downlabeledarrow[\big]
{\cs{free loops}\cs{classifying map}_{\#}}{}
&&\downlabeledarrow[\big]{\cs{classifying map}_{\#}}{}\\
\cs{Lifts}{X}{\cs{classifying map} \circ f}{\cs{classifying map range}}
{\cs{free loop projection}}{\cs{free loops}\cs{classifying map range}} &
\rightlabeledarrow{\ \iota^2\ }{}&
[X, \cs{free loops}\cs{classifying map range}]&
\rightlabeledarrow{\ \cs{free loop projection}^2_{\#}\ }{}&
[X, \cs{classifying map range}]
\\
\end{matrix}
\]
commutes.
Hence, if the group 
$\cs{Lifts}{X}{\cs{classifying map} \circ f}{\cs{classifying map range}}{\cs{free loop projection}}{\cs{free loops}\cs{classifying map range}}$
is abelian the cokernel of $\cs{lift homomorphism}{\cs{classifying map}}$ 
can be worked out from knowledge of just the right-hand square. 
Specifically 
\begin{Corollary}\namedLabel{free loops calculation}
With notation as above, suppose 
$\cs{Lifts}{X}{\cs{classifying map} \circ f}{\cs{classifying map range}}{\cs{free loop projection}}{\cs{free loops}\cs{classifying map range}}$
is abelian. 
The set $\mathcal C = 
({\cs{free loop projection}^2_{\#}})^{-1}(\cs{classifying map} \circ f) 
\subset [X,\cs{free loops}\cs{classifying map range}]$ is a group.
The set $\cs{free loops}{\cs{classifying map}_{\#}}
\bigl(({\cs{free loop projection}^1_{\#}})^{-1}(f)\bigr)$ is a subgroup of 
$\mathcal C$ and there is a bijection between the coset space of this inclusion and 
the cokernel of $\cs{lift homomorphism}{\cs{classifying map}}$. 
\end{Corollary}

\def\lrtcmrm{\mu_2}
\newCommand{classifying map range multiplication}{\lrtcmrm}
\def\lrtcmrlm{\mu_3}
\newCommand{classifying map range loop multiplication}{\lrtcmrlm}

\section{Some results on H-spaces}
To go further with the analysis in the last section requires some hypotheses. 
Let $\cs{classifying map domain}$ and $\cs{classifying map range}$ be H-spaces which have the homotopy type 
of CW complexes. 
Do not assume that the classifying map
$\cs{classifying map}\colon \cs{classifying map domain} \to 
\cs{classifying map range}$ is an H-map. 
\namedRef{isotropy subgroup} under these additional assumptions was obtained 
by J.~W.~Rutter \cite{Rutter}*{Theorem 1.3.1,p.~382} and 
there is considerable overlap between his \S1.4 and the material here. 

If $Y$ has the homotopy type of a CW complex, so do 
$\cs{free loops}{Y}$ and $\cs{based loop space}{Y}$, 
see Milnor \cite{Milnor}*{Theorem 3, p.~276}. 
If $Y$ is an H-space, the section map 
$\cs{free loops section}\colon Y \to \cs{free loops}{Y}$
and the inclusion map 
$\cs{loop inclusion}{}\colon\cs{based loop space}{Y} \to 
\cs{free loops}{Y}$ can be multiplied using the H-space product
to give homotopy equivalences,
$\cs{based loop space}{Y} \times Y \to\cs{free loops}{Y}$, see 
James \& Thomas, \cite{James-Thomas}*{Theorem 2.7, p.~494}, or   
Zabrodsky, \cite{Zabrodsky}*{1.3.6 Proposition, p.~24}. 
It follows that for any $h\colon X \to Y$, 
$\cs{Lifts}{X}{h}{Y}{\cs{free loop projection}}{\cs{free loops}Y}
$
is isomorphic as a group to $[X,\cs{based loop space}{Y}]$. 
Since $Y$ is an H-space, $[X,\cs{based loop space}{Y}]$ 
is abelian and 
$\cs{Lifts}{X}{h}{Y}{\cs{free loop projection}}{\cs{free loops}Y}
= [X,\cs{based loop space}{Y}] \times h \subset
[X,\cs{based loop space}{Y}]\times [X, Y] = [X,\cs{free loops}{Y}]$.

Hence it suffices to understand $\cs{free loops}{\cs{classifying map}}_{\#}$ 
for $\cs{classifying map} \colon \cs{classifying map domain} \to 
\cs{classifying map range}$. 
If $\alpha\in[X, \cs{based loop space}{ \cs{classifying map domain}}]$ and 
$\beta\in [X,  \cs{classifying map domain}]$ write
$\alpha\times\beta$ for $\mu_i\bigl(
(\cs{loop inclusion}{i})_{\#}(\alpha),
(\cs{free loops section}_i)_{\#}(\beta)
\bigr)$
Hence, 
to understand $\cs{free loops}{\cs{classifying map}}_{\#}$ it suffices 
to understand $\cs{free loops}{\cs{classifying map}}_{\#}(\alpha\times\beta)$
where $\alpha\in[X, \cs{based loop space}{\cs{classifying map domain}}]$ and $\beta\in [X,\cs{classifying map domain}]$.

Zabrodsky \cite{Zabrodsky}*{\S 1.4, p.~25} discusses the deviation from a map 
being an H-map.
In this case, the deviation is a map 
$D\colon \cs{free loops}{\cs{classifying map domain}}\wedge \cs{free loops}{\cs{classifying map domain}} 
\to \cs{free loops}{\cs{classifying map range}}$ which depends on $\cs{classifying map}$ and is 
null-homotopic if and only if $\cs{classifying map}$ is an H-map. 

With $\alpha\in[X, \cs{based loop space}{\cs{classifying map domain}}]$ and $\beta\in [X,\cs{classifying map domain}]$ 
define $W(\alpha,\beta)$ as the composition
\[
X \rightlabeledarrow{\Delta}{} X\wedge X 
\rightlabeledarrow{\alpha\wedge\beta}{} 
\cs{based loop space}{\cs{classifying map domain}}\wedge \cs{classifying map domain}
\rightlabeledarrow{\cs{loop inclusion}{1} \wedge \cs{free loop projection}_1}{}
\cs{free loops}{\cs{classifying map domain}}\wedge \cs{free loops}{\cs{classifying map domain}} 
\rightlabeledarrow{D}{} \cs{free loops}{\cs{classifying map range}}
\]
Then
$
\cs{classifying map range multiplication}\bigl(W(\alpha,\beta), 
\cs{free loops}{\cs{classifying map}}_{\#}(
\alpha \times \beta)\bigr)
= (\cs{based loop space}{\cs{classifying map}})_{\#}(\alpha) \times 
\cs{classifying map}_{\#}(\beta)$. 

Assume further that $\cs{classifying map range}$ is homotopy-associative so that 
$\cs{free loops}{\cs{classifying map range}}$ is also homotopy-associative. 
Then $[X,\cs{free loops}{\cs{classifying map range}}]$ is a group and so
\[\cs{free loops}{\cs{classifying map}}_{\#}(
\alpha \times \beta)\bigr) = 
\cs{classifying map range multiplication}\bigl(W(\alpha,\beta)^{-1}, 
(\cs{based loop space}{\cs{classifying map}})_{\#}(\alpha) \times 
\cs{classifying map}_{\#}(\beta)\bigr)~.\]

To continue, Zabrodsky \cite{Zabrodsky}*{1.4.2 Proposition, p.~25} shows
that \hskip 4pt
\[\begin{matrix}
\cs{free loops}{\cs{classifying map domain}} \wedge
\cs{free loops}{\cs{classifying map domain}} 
&\rightlabeledarrow{D}{}& \cs{free loops}{\cs{classifying map range}}\\
\downlabeledarrow{\cs{free loop projection}_1 \wedge \cs{free loop projection}_1}{}
&&\downlabeledarrow{\cs{free loop projection}_2}{}\\
\cs{classifying map domain}\wedge \cs{classifying map domain}&\rightlabeledarrow{D}{}&\cs{classifying map range}\\
\end{matrix}
\]
commutes. 
Hence it follows that the composition
$\cs{based loop space}{\cs{classifying map domain}} \wedge \cs{classifying map domain}
\rightlabeledarrow{\cs{loop inclusion}{1} \wedge \cs{free loop projection}_1}{}
\cs{free loops}{\cs{classifying map domain}}\wedge 
\cs{free loops}{\cs{classifying map domain}} 
\rightlabeledarrow{D}{} \cs{free loops}{\cs{classifying map range}}$ lifts to a map
$\cs{based loop space}{\cs{classifying map domain}} \wedge \cs{classifying map domain} 
\to \cs{based loop space}{\cs{classifying map range}}$.
This is a map into $\cs{based loop space}{\cs{classifying map range}}$ so it has a 
multiplicative inverse
$\mathfrak D\colon 
\cs{based loop space}{\cs{classifying map domain}}\wedge \cs{classifying map domain} \to \cs{based loop space}{\cs{classifying map range}}$. 
Further, for $\alpha\in[X, \cs{classifying map domain}]$ and $\beta\in [X,\cs{based loop space}{\cs{classifying map domain}}]$
define $\alpha\wedge_{\cs{classifying map}}\beta$ as the composition
$X \rightlabeledarrow{\Delta}{} X\wedge X 
\rightlabeledarrow{\alpha\wedge\beta}{} 
\cs{based loop space}{\cs{classifying map domain}}\wedge \cs{classifying map domain}
\rightlabeledarrow{\mathfrak D}{}
\cs{based loop space}{\cs{classifying map range}}$. 
Note $\alpha\wedge_{\cs{classifying map}}\beta$ is bilinear in both
$\alpha$ and $\beta$.

Plugging this into the formula above shows 
$\cs{free loops}{\cs{classifying map}}_{\#}(
\alpha \times \beta)\bigr) = 
\cs{classifying map range multiplication}\bigl(
(\cs{loop inclusion}{2})_{\#}(
\alpha\wedge_{\cs{classifying map}} \beta), 
(\cs{based loop space}{\cs{classifying map}})_{\#}(\alpha) \times 
\cs{classifying map}_{\#}(\beta)\bigr)$. 
Let $\cs{classifying map range loop multiplication}\colon 
\cs{based loop space}{\cs{classifying map range}} \times 
\cs{based loop space}{\cs{classifying map range}}
\to \cs{based loop space}{\cs{classifying map range}}$ be the usual H-space multiplication and since 
$\cs{classifying map range multiplication}$ is 
homotopy-associative the next formula has been proved:
\newNumber{formula for computing L_omega_sharp}
\[
\cs{free loops}{\cs{classifying map}}_{\#}(
\alpha \times \beta)\bigr) = 
\cs{classifying map range loop multiplication}
\bigl( (\alpha\wedge_{\cs{classifying map}}\beta), 
(\cs{based loop space}{\cs{classifying map}})_{\#}(\alpha)\bigr) \times 
\cs{classifying map}_{\#}(\beta)
\leqno(\ref{formula for computing L_omega_sharp})
\]
Formula \ref{formula for computing L_omega_sharp}, 
\namedRef{free loops calculation} 
and \namedRef{isotropy subgroup} prove
\begin{Theorem}\namedLabel{isotropy subgroup with H-spaces}
Let $\cs{classifying map domain}$ and $\cs{classifying map range}$ be H-spaces with $\cs{classifying map range}$ homotopy-associative.
Let $\cs{classifying map}\colon \cs{classifying map domain} \to \cs{classifying map range}$ be any map.
Let $\cs{classifying map total space}$ be the homotopy fibre of $\cs{classifying map}$, so
$\cs{based loop space}{\cs{classifying map range}} \to 
\cs{classifying map total space} \rightlabeledarrow{\cs{fibration projection}}{} \cs{classifying map domain}$ 
is a principal fibration. 
Let $\beta\in [X,\cs{classifying map domain}]$ be such that $\cs{classifying map}_{\#}(\beta) = 0$. 
Then $(\cs{fibration projection}_{\#})^{-1}(\beta) \subset 
[X,\cs{classifying map total space}]$ 
is non-empty and there is a bijection 
between $(\cs{fibration projection}_{\#})^{-1}(\beta)$ and the cokernel of the 
homomorphism
$\psi\colon [X,\cs{based loop space}{\cs{classifying map domain}}] \to [X,\cs{based loop space}{\cs{classifying map range}}]$ 
defined by 
$\psi(\alpha) = \cs{classifying map range loop multiplication}
\bigl((\alpha\wedge_{\cs{classifying map}}\beta), 
(\cs{based loop space}{\cs{classifying map}})_{\#}(\alpha)\bigr)$ for each 
$\alpha\in [X,\cs{based loop space}{\cs{classifying map range}}]$. 
\end{Theorem}

\begin{Remark}\namedLabel{free loops calculation II}
Continuing in this vein, let $e\in [X,\cs{classifying map total space}]$ 
be some element with
${\mathfrak p}_{\#}(e) = \beta$. 
Let 
\[\mathfrak D^\prime\colon 
\cs{based loop space}{\cs{classifying map domain}} \wedge \cs{classifying map total space} 
\rightlabeledarrow{\ 1_{\cs{based loop space}{\cs{classifying map domain}}}\wedge\, \mathfrak p\ }{}
\cs{based loop space}{\cs{classifying map domain}}\wedge \cs{classifying map domain}
\rightlabeledarrow{\mathfrak D}{} \cs{based loop space}{\cs{classifying map range}}
\]
and define $\alpha\wedge^\prime_{\cs{classifying map}} e$ as the composition
$X \rightlabeledarrow{\Delta}{} X\wedge X 
\rightlabeledarrow{\ \alpha\wedge e\ }{} 
\cs{based loop space}{\cs{classifying map domain}}\wedge \cs{classifying map total space}
\rightlabeledarrow{\mathfrak D^\prime}{}
\cs{based loop space}{\cs{classifying map range}}$. 
Certainly 
$\cs{classifying map range loop multiplication}
\bigl((\alpha\wedge_{\cs{classifying map}}\beta), 
(\cs{based loop space}{\cs{classifying map}})_{\#}(\alpha)\bigr)$
and 
$\cs{classifying map range loop multiplication}
\bigl((\alpha\wedge^\prime_{\cs{classifying map}}e), 
(\cs{based loop space}{\cs{classifying map}})_{\#}(\alpha)\bigr)
$ have the same image and sometimes $\mathfrak D^\prime$ 
is easier to compute than $\mathfrak D$.
\end{Remark}

Further information on $\mathfrak D$ 
can be obtained by applying 
(\ref{formula for computing L_omega_sharp}) to the identity map
which yields the next result.

\begin{Theorem}\namedLabel{determine D}
The composition 
$\cs{based loop space}{\cs{classifying map domain}}\times \cs{classifying map domain} 
\rightlabeledarrow{}{}
\cs{free loops}{\cs{classifying map domain}}
\rightlabeledarrow{\cs{free loops}{\cs{classifying map}}}{}
\cs{free loops}{\cs{classifying map range}}$ 
is homotopic to the following composition.
\[
\cs{based loop space}{\cs{classifying map domain}}\times \cs{classifying map domain} 
\rightlabeledarrow{\cs{smash quotient map}\times 1}{}
(\cs{based loop space}{\cs{classifying map domain}}\wedge \cs{classifying map domain})\times
(\cs{based loop space}{\cs{classifying map domain}}\times \cs{classifying map domain})
\rightlabeledarrow{\mathfrak D \times \cs{based loop space}{\cs{classifying map}}
\times \cs{classifying map}}{}
\cs{based loop space}{\cs{classifying map range}}\times \cs{based loop space}{\cs{classifying map range}}\times \cs{classifying map range}
\rightlabeledarrow{}{}
\cs{based loop space}{\cs{classifying map range}}\times \cs{classifying map range}
\rightlabeledarrow{}{}
\cs{free loops}{\cs{classifying map range}}
\]
where $\cs{smash quotient map}\colon 
\cs{based loop space}{\cs{classifying map domain}}\times \cs{classifying map domain} \to
\cs{based loop space}{\cs{classifying map domain}}\wedge \cs{classifying map domain}$ is the usual map. 
\end{Theorem}

\begin{Corollary}\namedLabel{determine D prime}
Suppose $a\in H_{r_1}(\cs{based loop space}{\cs{classifying map domain}};\Z)$ and
$b \in H_{r_2}(\cs{classifying map total space};\Z)$ are primitive classes, $r_i>0$. 
Then $\mathfrak D^\prime_\ast(a\times b) \in 
H_{r_1+r_2}(\cs{based loop space}{\cs{classifying map range}};\Z)$ maps to
$\cs{free loops}{\cs{classifying map}}_\ast(a\times b)\in 
H_{r_1+r_2}(\cs{free loops}{\cs{classifying map range}};\Z)$. 
\end{Corollary}
\begin{proof}
Since both $a$ and $b$ are primitive,  the composition 
\[\cs{based loop space}{\cs{classifying map domain}} \times \cs{classifying map total space}
\rightlabeledarrow{\Delta}{}
(\cs{based loop space}{\cs{classifying map domain}} \times \cs{classifying map total space}) 
\times 
(\cs{based loop space}{\cs{classifying map domain}} \times \cs{classifying map total space}) 
\to
(\cs{based loop space}{\cs{classifying map domain}} \wedge 
\cs{classifying map total space}) \times 
(\cs{based loop space}{\cs{classifying map domain}} \times \cs{classifying map total space})
\]
on $a\times b$ is $(a\wedge b)\times (1\times 1) + 1\times( a\times b)$. 
By \namedRef{determine D} the result follows since $\cs{classifying map}_\ast(b)=0$.
\end{proof}

\begin{Remark}\namedLabel{loop suspension remark}
If $\cs{classifying map range}$ is not an H-space but is highly connected,
then replace 
$\cs{classifying map range}$ by $\cs{loop suspension}{\cs{classifying map range}}$
and consider the composition
$\cs{classifying map domain}\to 
\cs{classifying map range} \rightlabeledarrow{\iota}{}
\cs{loop suspension}{\cs{classifying map range}}
$ where $\iota$ is the canonical inclusion. 
There is a commutative ladder
\[
\begin{matrix}
\cs{based loop space}{\cs{classifying map range}}& \to&
\cs{classifying map total space}&\to&\cs{classifying map domain}&
\rightlabeledarrow{\cs{classifying map}}{}&
\cs{classifying map range}\\
\downlabeledarrow{\cs{based loop space}{\iota}}{}&&
\downlabeledarrow{}{}&&\downlabeledarrow{}{}&&\downlabeledarrow{\iota}{}\\
\cs{loop suspension}[2]{\cs{classifying map range}}&\to&
\hat{\cs{classifying map total space}}&\to&\cs{classifying map domain}&
\rightlabeledarrow{}{}&
\cs{loop suspension}{\cs{classifying map range}}\\
\end{matrix}
\]
If $\pi_i(\cs{classifying map range}) = 0$ for $i<n$, then 
for any complex $X$ of dimension $\leqslant 2n-2$,
$[X,\cs{classifying map range}] \to [X,\cs{loop suspension}{\cs{classifying map range}}]$
is an isomorphism as are the other induced vertical maps. 
The results above can be applied to the $\hat{\cs{classifying map total space}}$ 
principal fibration to yield 
results about the $\cs{classifying map total space}$ principal fibration. 
\end{Remark}

\section{Some examples}
\subsection{Steenrod's problem}
Steenrod \cite{Steenrod} solved the problem of enumerating the homotopy
classes of maps $[X, S^n]$ where $n\geqslant 3$ and $X$ is a CW complex
of dimension at most $n+1$. 
\namedRef{isotropy subgroup} is not needed for the calculations in this subsection, but 
the results are needed below. 
A modern approach to this problem goes as follows. 

\def\SE#1{SE_{#1}}
For $n\geqslant 1$, let 
$\SE{n}$ be the fibre of the map $K(\Z,n) \rightlabeledarrow{Sq^2}{}
K(\cy{2},n+2)$. 
There is a map $S^n \to \SE{n}$ and the induced map
$[X, S^n] \to [X,\SE{n}]$ is an isomorphism if $n\geqslant 3$ and the dimension
of $X$ is at most $n+1$. 
In other words, $\SE{n}$ is the first two stages of a Postnikov decomposition for 
$S^n$.
The needed calculations are due to Serre \cite{Serre}. 

For $n\geqslant3$, $\SE{n} = \Omega\SE{n+1}$ so $\SE{n}$ is a homotopy-abelian
H-space, $[X,\SE{n}]$ is an abelian group, and the fibration de-loops. 
Write $\coker{\overline{Sq}^2}$ for the $\cy{2}$ vector space 
$H^{n+1}(X;\cy{2})/Sq^2\bigl(H^{n-1}(X;\Z)\bigr)$. 
Steenrod's main theorem 
\cite{Steenrod}*{Theorem 28.1, p.~318} follows:
\newNumber{Steenrod's answer}
\[0\to\coker{\overline{Sq}^2}
\to [X,S^{n}] \to H^n(X;\Z)\to 0\leqno(\ref{Steenrod's answer})\] 
is an exact sequence of abelian groups. 

Historically of course this approach is backwards. 
Steenrod invented $Sq^2$ to solve this problem and then worked out the 
Steenrod algebra which led to Serre's work. 
One could make a case for this being one of the all-time most important 
problems in algebraic topology. 

Larmore and Thomas \cite{Larmore-Thomas}*{\S5} give a procedure to
determine the extension in (\ref{Steenrod's answer}). 
In this case their procedure reduces to determining how the kernel of the 
multiplication by $2^k$ on $H^n(X; \Z)$ maps into $[X,S^n]$. 
To analyze this, consider the $2^k$ power maps on $\SE{n}$, $\cs{power map}k$,
$k\geqslant 1$. 
For each $k$ there is a commutative ladder of fibrations

\begin{hLadder}{1}3,2.5,3,\\
{K(\cy{2},n+1)}{\times 0}{K(\cy{2},n+1)}
{\SE{n}}{\cs{power map}k}{\SE{n}}
{K(\Z,n)}{\times 2^k}{K(\Z,n)}
{K(\cy{2},n+2)}{\times 0}{K(\cy{2},n+2)}
\\
{}{}
{}{}
{Sq^2}{Sq^2}
\\
\draw[dotted,->] (31) to node [above] {$\scriptstyle\cs{power lift}k$} (22); 
\end{hLadder}

Since the rows are fibrations (up to homotopy) there exists a map
$\cs{power lift}k$ as indicated in the diagram making the lower triangle
commute. 
Since $H^{n+1}\bigl(K(\Z,n);\cy{2^k}\bigr) = 0$, the map $\psi^\prime$ is 
unique.  
It follows from the Serre spectral sequence for the fibration that
$H^{n+1}\bigl(\SE{n};\cy{2^k}\bigr) = 0$ so the upper triangle involving 
$\cs{power lift}k$ also commutes. 

Next check that the following diagram commutes. 

\hskip 1in\begin{vLadder}{4}1.1,1.2,\\
{K(\Z,n)}{\cs{power lift}k}{\SE{n}}
{K(\Z,n)}{\times 1}{K(\Z,n)}
{K(\cy{2^k},n)}{Sq^2}{K(\cy{2},n+2)}
\\
{\times 2^k}{}
{}{Sq^2}
\\
\end{vLadder}

It follows that there is an induced map on the fibres 
which is the loops of $Sq^2$ and is therefore again $Sq^2$. 
Hence

\hskip 1in\begin{vLadder}{4}1.1,1.2,\\
{K(\cy{2^k},n-1)}{Sq^2}{K(\cy{2},n+1)}
{K(\Z,n)}{\cs{power lift}k}{\SE{n}}
{K(\Z,n)}{\times 1}{K(\Z,n)}
\\
{\delta_k}{}
{2^k}{}
\\
\end{vLadder}

\noindent
commutes, where $\delta_k$ is the evident Bockstein. 
The next result summarizes the above discussion.

\def\ainthm{\gamma}
\begin{Theorem}\namedLabel{determine the extension}
Let $X$ be a finite complex of dimension $\leqslant n+1$. 
Fix $\ainthm\in H^n(X;\Z)$ and suppose there is a $k\geqslant 1$ 
such that $2^k \ainthm = 0$.  
Pick $\ainthm ^\prime\in H^{n-1}(X;\cy{2^k})$
with $\delta_k(\ainthm ^\prime) = \ainthm $ and then compute
$Sq^2(\ainthm ^\prime)\in H^{n+1}(X;\cy{2})/Sq^2\bigl(H^{n-1}(X;\Z)\bigr) 
\subset [X,S^n]$.
For any $\bar{\ainthm}\in [X,S^n]$ which maps to $\ainthm $,
$2^k \bar{\ainthm} = \cs{power lift}k(\ainthm) = Sq^2(\ainthm ^\prime)$.
\end{Theorem}

\begin{Example}
Suppose $X$ is a complex of dimension $\leqslant n+1$ and  suppose that
$Sq^2\colon H^{n-1}(X;\Z) \to H^{n+1}(X;\cy{2})$ and
$Sq^2\colon H^{n-1}(X;\cy{2}) \to H^{n+1}(X;\cy{2})$ have the same image.
Then 
$[X,S^n] = \coker{\overline{Sq}^2}\oplus H^n(X;\Z)$. 
\end{Example}

\begin{Example}
If $X^4$ is Habegger's manifold \cite{Habegger} or an Enriqu\'e's surface, 
then $Sq^2\colon H^2(X;\Z) \to H^4(X;\cy{2})$ is
zero but  $Sq^2\colon H^2(X;\cy{2}) \to H^4(X;\cy{2})$ is onto.
Since $H^3(X;\Z) = \cy{2}$ it follows that $[X,S^3] \cong \cy{4}$. 
\end{Example}

\subsection{Pontrjagin's problem}
Pontrjagin \cite{Pontrjagin} solved the problem of enumerating $[X,S^2]$ 
for $X$ a 3-complex before Steenrod did his work. 
From the point of view taken here, $S^2\to \cs{classifying space}{S^1} 
\to \cs{classifying space}{S^3}$ 
is a fibration so $S^2$ is the total space of a principal fibration,
$S^3 \to S^2 \to \cs{classifying space}{S^1}$. 
Since $S^1$ is an abelian group, $\cs{classifying space}{S^1} = \cs{CP}{\infty}$ is an
H-space.
However, $S^3$ is not abelian and $\cs{classifying space}{S^3} = \cs{HP}{\infty}$ is not
an H-space. 

However, $\pi_i(\cs{classifying space}{S^3}) = 0$ for $i< 4$ so 
\namedRef{loop suspension remark} says that as long as the 
dimension of $X$ is $\leqslant 2\cdot 4 -2 = 6$, the theorems in \S 5 apply.
The next subsection computes the answer for all complexes of 
dimension $\leqslant 4$ and includes a statement and proof of Pontrjagin's result
as \namedRef{Pontrjagin's calculation}. 

\subsection{The second cohomotopy set of a 4-complex}\label{second cohomotopy set}
Let $X$ have the homotopy type of a CW-complex of dimension $\leqslant 4$. 
The first step is to compute the map 
$[X, \cs{classifying space}{S^1}] \to [X,\cs{classifying space}{S^3}]$.
The map $\cs{classifying space}{S^3} \to K(\Z,4)$ 
giving a generator of $H^4(\cs{classifying space}{S^3};\Z)\cong\Z$
is 5-connected, so $[X,\cs{classifying space}{S^3}] \to [X,K(\Z,4)] = H^4(X;\Z)$ 
is an isomorphism.
Since the map
$\cs{classifying space}{S^1} \to \cs{classifying space}{S^3}$ 
is the standard inclusion of $\cs{CP}{\infty}$ in
$\cs{HP}{\infty}$, the map 
$[X,\cs{classifying space}{S^1}] = H^2(X;\Z) \to 
[X,\cs{classifying space}{S^3}]=H^4(X;\Z)$ 
is just the cup product square. 
Hence $[X,S^2] \to H^2(X;\Z)$ is onto the subset of classes $\beta\in H^2(X;\Z)$ 
such that $\beta\cup \beta = 0\in H^4(X;\Z)$. 

Since $\cs{classifying space}{S^3}$ is not an H-space, use 
\namedRef{loop suspension remark} and work with
$\cs{loop suspension}{\cs{classifying space}{S^3}}$.

In \S6.1, the group $[X,S^3] = [X,\cs{loop suspension}[2]{\cs{classifying space}{S^3}}]$
was computed for any 4-complex. 

\medskip

For a fixed map $e\colon X \to S^2$, the next step is to understand 
the homomorphism
$\cs{homomorphism to [X,S^3]}[e]\colon H^1(X;\Z) \to 
[X,\cs{loop suspension}{\cs{classifying space}{S^3}}]$. 
Since $\alpha\in H^1(X;\Z)$ is equivalent to a homotopy class of based maps
$\alpha\colon X \to S^1$,
and since 
$\mathfrak D^\prime \colon S^1\wedge S^2 \to 
\cs{loop suspension}[2]{\cs{classifying space}{S^3}}$, it follows that  
$\mathfrak D^\prime$ factors through the degree $c_{e}$-map
$S^3 \to S^3$.
Hence there is a homomorphism 
$\bar{\psi}\colon H^1(X;\Z) \to [X,S^3]$ such that
$\cs{homomorphism to [X,S^3]}[e]$ is the composition
\[
H^1(X;\Z) \rightlabeledarrow{\bar{\psi}}{} [X,S^3] 
\rightlabeledarrow{({c_{e}})_{\#}}{} [X, S^3]\] 
where $(c_{e})_{\#}$ is the map induced by the 
degree $c_e$ map on $S^3$. 
Since $[X,S^3]$ is an abelian group, $({c_{e}})_{\#}$ is just multiplication 
by $c_{e}$.

By definition, the composition $H^1(X;\Z) \rightlabeledarrow{\bar{\psi}}{}
[X,S^3] \to H^3(X;\Z)$ just sends $\alpha$ to $\alpha \cup \beta$ where
$\beta\in H^2(X;\Z)$ is given by pulling back the fundamental class in 
$H^2(S^2;\Z)$ via $e\colon X \to S^2$. 
It follows from \namedRef{free loops calculation lemma} below that $c_{e} = \pm 2$. 
The sign will not be determined here. 
\begin{Lemma}\namedLabel{free loops calculation lemma}
The map 
$H_3(\cs{free loops}{\cs{classifying space}{S^1}};\Z) \to 
H_3(\cs{free loops}{\cs{classifying space}{S^3}};\Z)$ is 
multiplication by $\pm 2$.
\end{Lemma}
\begin{proof}
For $m=1$ or $3$, the Serre spectral sequence for $S^m \to 
\cs{free loops}{\cs{classifying space}{S^m}} \to
\cs{classifying space}{S^m}$ collapses and 
$H_\ast(\cs{free loops}{\cs{classifying space}{S^m}};\Z) = E(e_m)\otimes \Z[x_{m+1}]$  
where $e_m\in H_m(\cs{free loops}{\cs{classifying space}{S^m}};\Z)$ is the image of 
$H_m(S^m;\Z)$; $x_{m+1}\in H_{m+1}(\cs{free loops}
{\cs{classifying space}{S^m}};\Z)$ maps 
to a generator of $H_{m+1}(\cs{classifying space}{S^m};\Z)$; $E(e_m)$ 
is an exterior algebra
and $\Z[x_{m+1}]$ is a polynomial algebra. 

Now $H_{2}(\cs{free loops}{S^2};\Z) \cong \Z\oplus\cy{2}$, 
say by Ziller, \cite{Ziller}*{p.~21}. 
It follows that in the Serre spectral sequence for the fibration 
$\cs{free loops}{S^2} \to \cs{free loops}{\cs{classifying space}{S^1}} \to 
\cs{free loops}{\cs{classifying space}{S^3}}$ 
there is a single differential from 
$H_3(\cs{free loops}{\cs{classifying space}{S^3}};\Z)$ onto $\cy{2}$ so
$H_3(\cs{free loops}{\cs{classifying space}{S^1}};\Z) \to 
H_3(\cs{free loops}{\cs{classifying space}{S^3}};\Z)$ is 
multiplication by $\pm2$. 
\end{proof}

It follows that the homomorphism 
$\cs{homomorphism to [X,S^3]}[e]\colon H^1(X;\Z) \to [X,S^3]$ factors as 
\[
H^1(X;\Z) \rightlabeledarrow{\ (\, \underline{\hskip 6pt}\, )\, \cup\,\beta\ }{} 
H^3(X;\Z) \rightlabeledarrow{\ \cs{power lift}1\ }{} [X,S^3]\]  
and so
$\cs{homomorphism to [X,S^3]}[e]$ only depends on $\beta$ and hereafter will
be written $\cs{homomorphism to [X,S^3]}[\beta]$
\begin{Theorem}\namedLabel{the second cohomotopy set}
Let $X$ be a complex of dimension $\leqslant 4$ and let 
$p_{\#}\colon [X,S^2] \to H^2(X;\Z)$ be the map pulling back a fixed
generator of $H^2(S^2;\Z)$. 

If $\beta\in H^2(X;\Z)$ is given, then $p_{\#}^{-1}(\beta)$ is non-empty 
if and only if $\beta \cup \beta = 0$. 
Furthermore, if $p_{\#}^{-1}(\beta)$ is non-empty, then there is a bijection between it
and the cokernel of $\cs{homomorphism to [X,S^3]}[\beta]\colon 
H^1(X;\Z) \to [X,S^3]$. 
\end{Theorem}

\begin{Remark}
Let $\cs{cohomology cokernel}{\beta}$ be the cokernel of
$H^1(X;\Z) \rightlabeledarrow{2(\ \underline{\hskip 5pt}\ )\,\cup\,\beta\ }{}H^3(X;\Z)$. 
Then there is an exact sequence 
\[
\coker{\overline{Sq}^2} \rightlabeledarrow{\ q\ }{} 
\coker{\cs{homomorphism to [X,S^3]}[\beta]} 
\to \cs{cohomology cokernel}{\beta}\to0\]
The kernel of $q$ is the set of all elements of the form $Sq^2(a)$ 
for some $a\in H^2(X;\cy{2})$ such that there exists $\alpha\in H^1(X;\Z)$
such that $\delta_1(a) = \alpha\cup\beta \in H^3(X;\Z)$. 
\end{Remark}

\begin{Remark}
There are three types of connected, closed, compact 4-manifolds:
(1) there exists an $x\in H^2(X;\Z)$ with odd square;
(2) for all $x\in H^2(X;\cy{2})$ $x\cup x =0$;
(3) $X$ is not of type (2) but for all $x\in H^2(X;\Z)$ $x\cup x$ is even. 
If $X$ has type (1), $\coker{\cs{homomorphism to [X,S^3]}[\beta]} 
\to \cs{cohomology cokernel}{\beta}$ is an isomorphism. 
If $X$ has type (2) $0\to \cy{2} \to \coker{\cs{homomorphism to [X,S^3]}[\beta]} 
\to \cs{cohomology cokernel}{\beta}\to 0$ is split exact. 
If $X$ has type (3) $\cy{2} \to \coker{\cs{homomorphism to [X,S^3]}[\beta]} 
\to \cs{cohomology cokernel}{\beta}\to 0$ is exact and
\namedRef{determine the extension} can be used to determine the group.
If $X$ has type (3) and if $\coker{\cs{homomorphism to [X,S^3]}[\beta]} 
\to \cs{cohomology cokernel}{\beta}$ is not an isomorphism, then the 
sequence is not split.
The manifold $\cs{CP}2$ has type (1), any Spin manifold has type (2) 
and the Habegger manifold \cite{Habegger} is an example
a type (3) manifold for which the extension is not split. 
The author does not know an example of a type (3) manifold for which 
$\coker{\cs{homomorphism to [X,S^3]}[\beta]} 
\to \cs{cohomology cokernel}{\beta}$ is an isomorphism. 
\end{Remark}

\begin{Corollary}[Pontrjagin, \cite{Pontrjagin}]\namedLabel{Pontrjagin's calculation}
If $X$ is a 3-dimensional complex then $[X,S^2] \to H^2(X;\Z)$ is onto and 
there is a bijection between $p_{\#}^{-1}(\beta)$ and 
$\cs{cohomology cokernel}{\beta}$. 
\end{Corollary}

\begin{Example}\namedLabel{4-dimensional example}
Let $X = S^2 \times S^1$. 
Then $H^2(X;\Z) \cong \Z$: let $\gamma$ be a generator. 
If $\beta = c \gamma$ then there are maps $X \to S^2$ 
such that $\beta$ is the image of a generator of $H^2(S^2;\Z)$ 
and there is a bijection between $p_{\#}^{-1}(\beta)$ and $\Z$ if $c=0$ and 
$\cy{2 c}$ otherwise. 
\end{Example}

\begin{Example}\namedLabel{S^2 x torus}
Let $X = S^2 \times S^1\times S^1$. 
Let $\{\mathfrak a_1, \mathfrak a_2\} \subset H^1(X;\Z) \cong \Z\oplus\Z$ 
be a basis and let
$\{\mathfrak a = \mathfrak a_1\cup \mathfrak a_2, \mathfrak b\}\subset H^2(X;\Z)$ 
be a basis. 
It follows that $\{\mathfrak b \cup \mathfrak a_1, \mathfrak b \cup \mathfrak a_2\}$ 
is a basis for $H^3(X;\Z)$. 
Then $\beta = a \mathfrak a + b \mathfrak b$ has square $0$ 
if and only if $a\cdot b = 0$.
If $b = 0$, then $\coker{\cs{homomorphism to [X,S^3]}[\beta]} = 
H^3(X;\Z) \oplus \cy{2} \cong \Z^2\oplus \cy{2}$. 
If $a = 0$, then the image of $\cs{homomorphism to [X,S^3]}[\beta]$ 
is spanned by $(2 b)\, \mathfrak b \cup \mathfrak a_1$ and 
$(2 b)\, \mathfrak b \cup \mathfrak a_2$ and so 
$\coker{\cs{homomorphism to [X,S^3]}[\beta]}  \cong \cy{2 b} \oplus\cy{2 b} \oplus \cy{2}$.
\end{Example}


\begin{thebibliography}{WW} 
\bib{Habegger}{article}{
   author={Habegger, Nathan},
   title={Une vari\'et\'e de dimension $4$ avec forme d'intersection paire
   et signature $-8$},
   language={French},
   journal={Comment. Math. Helv.},
   volume={57},
   date={1982},
   number={1},
   pages={22--24},
   issn={0010-2571},
   review={\MR{672843 (83k:57018)}},
}
\bib{James-Thomas}{article}{
  author={James, Ioan},
  author={Thomas, Emery},
  title={An approach to the enumeration problem for non-stable vector
   bundles},
  journal={J. Math. Mech.},
  volume={14},
  date={1965},
  pages={485--506},
  review={\MR{0175134 (30 \#5319)}},
  }

\bib{Larmore-Thomas}{article}{
  author={Larmore, Lawrence L.},
  author={Thomas, Emery},
  title={Group extensions and principal fibrations},
  journal={Math. Scand.},
  volume={30},
  date={1972},
  pages={227--248},
  issn={0025-5521},
  review={\MR{0328935 (48 \#7277)}},
  }

\bib{Milnor}{article}{
  author={Milnor, John},
  title={On spaces having the homotopy type of a ${\rm CW}$-complex},
  journal={Trans. Amer. Math. Soc.},
  volume={90},
  date={1959},
  pages={272--280},
  issn={0002-9947},
  review={\MR{0100267 (20 \#6700)}},
  }

\bib{Nomura}{article}{
  author={Nomura, Yasutoshi},
  title={On mapping sequences},
  journal={Nagoya Math. J.},
  volume={17},
  date={1960},
  pages={111--145},
  issn={0027-7630},
  review={\MR{0132545 (24 \#A2385)}},
  }

\bib{Peterson}{article}{
   author={Peterson, Franklin P.},
   title={Functional cohomology operations},
   journal={Trans. Amer. Math. Soc.},
   volume={86},
   date={1957},
   pages={197--211},
   issn={0002-9947},
   review={\MR{0105679 (21 \#4417)}},
}

\bib{Peterson-Thomas}{article}{
   author={Peterson, Franklin P.},
   author={Thomas, Emery},
   title={A note on non-stable cohomology operations},
   journal={Bol. Soc. Mat. Mexicana (2)},
   volume={3},
   date={1958},
   pages={13--18},
   review={\MR{0105680 (21 \#4418)}},
}

\bib{Pontrjagin}{article}{
   author={Pontrjagin, Lev Semenovich},
   title={A classification of mappings of the three-dimensional complex into
   the two-dimensional sphere},
   language={English, with Russian summary},
   journal={Rec. Math. [Mat. Sbornik] N. S.},
   volume={9 (51)},
   date={1941},
   pages={331--363},
   review={\MR{0004780 (3,60g)}},
}

\bib{Rutter}{article}{
   author={Rutter, John W.},
   title={A homotopy classification of maps into an induced fibre space},
   journal={Topology},
   volume={6},
   date={1967},
   pages={379--403},
   issn={0040-9383},
   review={\MR{0214070 (35 \#4922)}},
}
\bib{Serre}{article}{
   author={Serre, Jean-Pierre},
   title={Cohomologie modulo $2$ des complexes d'Eilenberg-MacLane},
   language={French},
   journal={Comment. Math. Helv.},
   volume={27},
   date={1953},
   pages={198--232},
   issn={0010-2571},
   review={\MR{0060234 (15,643c)}},
}

\bib{Steenrod}{article}{
  author={Steenrod, Norman E.},
  title={Products of cocycles and extensions of mappings},
  journal={Ann. of Math. (2)},
  volume={48},
  date={1947},
  pages={290--320},
  issn={0003-486X},
  review={\MR{0022071 (9,154a)}},
  }

\bib{V}{article}{
  author={Vogt, Rainer M.},
  title={Convenient categories of topological spaces for homotopy theory},
  journal={Arch. Math. (Basel)},
  volume={22},
  date={1971},
  pages={545--555},
  issn={0003-889X},
  review={\MR{0300277 (45 \#9323)}},
  }
 
 \bib{Zabrodsky}{book}{
   author={Zabrodsky, Alexander},
   title={Hopf spaces},
   note={North-Holland Mathematics Studies, Vol. 22;
   Notas de Matem\'atica, No. 59},
   publisher={North-Holland Publishing Co.},
   place={Amsterdam},
   date={1976},
   pages={x+223},
   review={\MR{0440542 (55 \#13416)}},
}

\bib{Ziller}{article}{
  author={Ziller, Wolfgang},
  title={The free loop space of globally symmetric spaces},
  journal={Invent. Math.},
  volume={41},
  date={1977},
  number={1},
  pages={1--22},
  issn={0020-9910},
  review={\MR{0649625 (58 \#31198)}},
  }

\end{thebibliography}
\end{document}